\documentclass[12pt,reqno]{amsart}
\usepackage{amsmath,amssymb,amsfonts,amscd,latexsym,amsthm,mathrsfs,bigints,mathtools,verbatim,cite}
\usepackage[usenames]{color}
\usepackage{tikz,graphicx}
\usepackage{hyperref}
\newcommand{\wt}{\widetilde}
\renewcommand{\d}{{\mathrm d}}
\begin{document}

\title[Hypergeometric heritage of W.\,N.~Bailey]{Hypergeometric heritage of W.\,N.~Bailey.\\
With an appendix: Bailey's letters to F.~Dyson}

\author{Wadim Zudilin}
\address{IMAPP, Radboud Universiteit, PO Box 9010, 6500 GL Nijmegen, The Netherlands}
\email{w.zudilin@math.ru.nl}
\address{MAPS, The University of Newcastle, Callaghan, NSW 2308, Australia}
\email{wadim.zudilin@newcastle.edu.au}

\dedicatory{On the occasion of Bailey's 123rd birthday}


\begin{abstract}
We review some of W.\,N.~Bailey's work on hypergeometric functions that found solid applications in number theory.
The text is complemented by Bailey's letters to Freeman Dyson from the 1940s.
\end{abstract}

\maketitle

\section{Introduction}
\label{intro}

This paper is a review of those parts of the mathematical legacy of the late Wilfrid Norman Bailey (1893--1961) that have had greatest impact on
certain developments in number theory in which we have been involved. In other words, these are our personal encounters with the mathematics of Bailey,
mathematics\,---\,we believe\,---\,that should be better known. The details of Bailey's biography can be found in his obituary \cite{Sl62}
written by Lucy Slater, where she also summarises his research as follows:
\begin{quote}
\small
``His life work was in the field of classical analysis. He studied hypergeometric
functions, all the functions of mathematical physics which are
specialized cases of hypergeometric functions and the interconnections
between the various types of function. He was always seeking the underlying
common features of these functions and he had the aim to unify the
various theorems on special functions as far as possible into one theory of
general application.

He wrote one book, the Cambridge Tract \emph{Generalized Hypergeometric
Series}, by which his name will always be remembered. This was the first
work in English to be written on the subject. It gathered together all his
researches up to that date (1935) in a very readable form, together with
outlines of all the earlier work which had been carried out by his predecessors,
Gauss, Heine, Appell, and others, on the Continent. It has
been remarked by students that he did not always say who had proved
some of the theorems in the Tract, but, characteristically modest, he
replied: ``I did not want to keep on repeating the word Bailey.''
\end{quote}

The generalized hypergeometric series is defined as
\begin{equation}
{}_mF_{m-1}\biggl(\begin{matrix} a_1, \, a_2, \, \dots, \, a_m \\ b_2, \, \dots, \, b_m \end{matrix} \biggm|z\biggr)
=\sum_{n=0}^\infty\frac{(a_1)_n(a_2)_n\dotsb(a_m)_n}{(b_2)_n\dotsb(b_m)_n}\,\frac{z^n}{n!},
\label{eq01}
\end{equation}
where
$$
(a)_n=\frac{\Gamma(a+n)}{\Gamma(a)}=\begin{cases}
a(a+1)\dotsb(a+n-1) &\text{if $n\ge1$}, \\
1 &\text{if $n=0$},
\end{cases}
$$
denotes the Pochhammer symbol (or rising factorial). The series in \eqref{eq01} can be shown to be convergent inside the unit disk $|z|<1$,
and in that region defines an analytic function of $z$ that satisfies a (linear homogeneous) hypergeometric differential equation of order~$m$
with regualr singularities at $z=0$, $1$ and $\infty$. It is this differential equation which is commonly used to analytically continue the ${}_mF_{m-1}$ function defined in \eqref{eq01}
to the whole $\mathbb C$-plane with the cut between two singularities (for example, with the cut along $[1,\infty)$); the result
of such analytic continuation is the generalized ${}_mF_{m-1}$-hypergeometric function.

Writing the right-hand side in \eqref{eq01} as the sum $\sum_{n=0}^\infty u(n)$ we can notice that $u(0)=1$ and
$$
\frac{u(n+1)}{u(n)}=\frac{(a_1+n)(a_2+n)\dotsb(a_m+n)}{(1+n)(b_2+n)\dotsb(b_m+n)}\,z,
$$
so that a way to characterise a hypergeometric series (sum or function) is to say that the quotient of two successive terms is a rational function of the index.
Of course, in our particular setting \eqref{eq01} we additionally impose the condition on the degrees of the numerator and denominator of the rational function to coincide
as well as the specific factor $1+n$ to be presented in the denominator. But any of these can be relaxed; for example, having more generally the factor $b_1+n$
in place of $1+n$ leads to what is called a bilateral hypergeometric function. One can further generalise such generalized hypergeometric functions to
generalized hypergeometric integrals, commonly referred to as Barnes-type integrals (because of presence of many other integral representations for hypergeometric
functions), and to basic hypergeometric functions when the quotient $u(n+1)/u(n)$ of two successive terms is a rational function of $q^n$, where the base $q$
becomes another parameter. One more step of generalization of that process is considering double and multiple hypergeometric functions.
We will witness many of such hypergeometric creatures later on.

The fact that the first great hypergeometric players, like Euler, Gauss and Riemann, also happened to be number theorists
hints at possible links between the two mathematics areas. However, a division of certain topics within mathematics always has been subjective:
for example, the theta functions of Jacobi were, for a long part, a topic of analysis, whereas these days many of us place them in number theory
or even algebraic geometry. As Slater's quotation reveals, Bailey's research was in analysis and his treatment of the functions he investigated was entirely analytical.
Although he did not care about number theory much, we hope he would be very pleased to learn what a great impact on different topics of number theory his
results have had.

The rest of this paper is organised as follows. In Section~\ref{zeta} we outline some links between Bailey's transformations of hypergeometric series
and the irrationality and linear independence results for the values of Riemann's zeta function $\zeta(s)=\sum_{n=1}^\infty n^{-s}$ at positive integers. Bailey's investigations
on reduction of Appell's (double) hypergeometric functions to single hypergeometrics are reviewed in Section~\ref{sun}; some number-theoretical consequences
of these results are new Ramanujan-type identities for $1/\pi$ which were recently discovered experimentally by the Chinese mathematician Zhi-Wei Sun.
In Section~\ref{qsect} we examine one of first signs of the `$q$-desease' in Bailey\,---\,a proof of what can regarded as a $q$-analogue
of Euler's identity $\zeta(2,1)=\zeta(3)$, where
$$
\zeta(2,1)=\sum_{n>l\ge1}\frac1{n^2l}.
$$
Our final Section~\ref{erdos} outlines some further connections of Bailey's research with number theory, in particular, his work on generalized
Rogers--Ramanujan identities with influences from and to Freeman Dyson, the work that eventually led to what George Andrews called the Bailey pairs and Bailey chains in 1984.
Apart from Bailey's famous and highly cited book \cite{B35}, the latter objects are the known hypergeometric heritage of W.\,N.~Bailey,
something that Slater could not forecast in 1962.

\medskip
\noindent
\textbf{Acknowledgements.}
Several people helped me to understand the personality of W.\,N.~Bailey, his time, and different levels of appreciation of hypergeometric functions
at different periods. It is my pleasure to thank George Andrews, Bruce Berndt, Frits Beukers, Gert Heckman, Freeman Dyson, Christian Krattenthaler, Jonathan Sondow, Ole Warnaar and one of the two anonymous referees of \emph{Bull.\ London Math.\ Soc.}\
for their valuable input and related comments on parts of the text. I am particularly thankful to Freeman Dyson for giving me access to his correspondence
with Bailey and permitting to include it as an appendix to this review.

\section{Transformations of hypergeometric series\\ and irrationality of zeta values}
\label{zeta}

It is hard to imagine that one innocent-looking identity \cite[eq.~(3.4)]{B32}, namely,
\begin{align}
&
{}_7F_6\biggl(\begin{alignedat}{7} & a, &\,& 1+\tfrac12a, &\,& \quad\; b, &\,& \qquad c, &\,& \qquad d, &\,& \qquad e, &\,& \qquad f \\[-3\jot]
 & &\,& \quad \tfrac12a, &\,& {\kern-2.5mm}1+a-b, &\,& 1+a-c, &\,& 1+a-d, &\,& 1+a-e, &\,& 1+a-f \end{alignedat} \biggm| 1 \biggr)
\nonumber\\ &\;
=\frac{\Gamma(1+a-b)\,\Gamma(1+a-c)\,\Gamma(1+a-d)\,\Gamma(1+a-e)\,\Gamma(1+a-f)}
{\splitfrac{\Gamma(1+a)\,\Gamma(b)\,\Gamma(c)\,\Gamma(d)\,\Gamma(1+a-c-d)\,\Gamma(1+a-b-d)}{\times\Gamma(1+a-b-c)\,\Gamma(1+a-e-f)}}
\nonumber\\ &\;\quad\times
\frac1{2\pi i}\,\raisebox{1.8mm}{$\displaystyle\bigintss_{-i\infty}^{i\infty}$}
\frac{\splitfrac{\Gamma(b+t)\,\Gamma(c+t)\,\Gamma(d+t)\,\Gamma(1+a-e-f+t)}{\times\Gamma(1+a-b-c-d-t)\,\Gamma(-t)\,\d t}}
{\Gamma(1+a-e+t)\,\Gamma(1+a-f+t)},
\label{7F6tr}
\end{align}
may have anything in common with the irrationality of $\zeta(3)$. And with much more irrational than that.

The hypergeometric series on the left-hand side in \eqref{7F6tr} is evaluated at $z=1$ and is of a very special type when its parameters come in pairs with the same sum:
$$
a+1=\bigl(1+\tfrac12a\bigr)+\tfrac12=b+(1+a-b)=c+(1+a-c)=\dots=f+(1+a-f).
$$
This property characterises a well-poised hypergeometric series\,---\,a term coined by F.\,J.\,W.~Whipple~\cite{Wh25}; the special form of the second
pair of parameter, $1+\tfrac12a$ and $\tfrac12a$, assigns it to the class of very-well-poised hypergeometric functions in modern terminology.
The majority of known summation and transformation formulae are specific to these classes, as they possess some additional structural symmetries.
The integral on the right-hand side in \eqref{7F6tr} is of Barnes type; the path of integration in it is parallel to the imaginary axis, except
that it is indented, if necessary, so that the decreasing sequences of poles of the integrand lie to the left, and the increasing sequences to the right of the contour.
If this hypergeometric integral is evaluated by considering the residues at poles on the right of the contour, then we obtain the transformation
of a very-well-poised ${}_7F_6$ in terms of two balanced (or Saalsch\"utzian) ${}_4F_3$; for a hypergeometric series the latter means that the sum
of the denominator parameters exceeds the sum of the numerator parameters by~1.

Before going into details of the interplay between \eqref{7F6tr} and the irrationality questions of the zeta values\,---\,the values
of Riemann's zeta function $\zeta(s)$ at integers $s=2,3,4,\dots$, let us make some historical remarks on the latter.
Giving the closed form $\zeta(2)=\pi^2/6$ for the first convergent zeta value was a part of Euler's triumph in his resolution of the Basel problem:
he also extended this more generally to $\zeta(2k)/\pi^{2k}\in\mathbb Q$ and, on this way, discovered the functional equation
that related $\zeta(s)$ to $\zeta(1-s)$. In absence of complex analysis at that time, Euler's interpretation of the divergent
series for $\zeta(s)$ when $s<0$ looks very impressive. It was much later, when functions of a complex variable became available,
that Euler's ideas and methods were put on a solid ground. That development was an essential
tool in Lindemann's 1882 proof \cite{Li82} that $\pi$ is a transcendental number. The latter result, together with Euler's,
left over a little mystery about the arithmetic nature of $\zeta(2),\zeta(4),\zeta(6),\dots$\,---\,the even zeta values\,---\,but
at the same time initiated an intrigue for the odd zeta values $\zeta(3),\zeta(5),\zeta(7),\dots$\,.
It was only in the 1970s when the senior mathematician Roger Ap\'ery from Caen, in a process of executing his convergence acceleration method
against a ``number table due to Ramanujan'' \cite{Ap96,Ra12}, managed to demonstrate that $\zeta(3)$ is irrational.
The detailed account of the controversial story of Ap\'ery's discovery \cite{Ap79} is given in the excellent
exposition \cite{vP79} by Alf van der Poorten. It took another couple of decades and the junior mathematician Tanguy Rivoal,
also from Caen, to prove in 2000 that infinitely many odd zeta values are irrational~\cite{Ri00}. The development
of the story from 1979 to 2000 was quite intense, with several interesting novelties appearing in both sharpening
of the number-theoretical tools and constructing rational approximations to zeta values. As we will see now,
the latter are all about hypergeometric functions.

Several ways are now known to write Ap\'ery's rational approximations $v_n/u_n\in\mathbb Q$, where $n=0,1,2,\dots$,
to~$\zeta(3)$. One can use the explicit formulae
\begin{equation}
u_n=\sum_{k=0}^n\binom nk^2\binom{n+k}k^2
\label{eq:0bin}
\end{equation}
and some more involved double sums for $v_n$ originally produced by Ap\'ery himself, or define sequences of both the denominators $\{u_n\}=\{u_n\}_{n=0,1,\dots}$
and of the numerators $\{v_n\}=\{v_n\}_{n=0,1,\dots}$ as solutions of the same polynomial recursion
\begin{subequations}
\label{AR}
\begin{equation}
(n+1)^3u_{n+1}-(2n+1)(17n^2+17n+5)u_n+n^3u_{n-1}=0
\label{eq:0.3}
\end{equation}
with the initial data
\begin{equation}
u_0=1, \quad u_1=5 \qquad\text{and}\qquad v_0=0, \quad v_1=6.
\label{eq:0.4}
\end{equation}
\end{subequations}
Alternatively, one can use the Beukers triple integral~\cite{Be79}
\begin{equation}
u_n\zeta(3)-v_n=\frac12\iiint\limits_{[0,1]^3}
\frac{x^n(1-x)^ny^n(1-y)^nz^n(1-z)^n}{(1-(1-xy)z)^{n+1}}\,\d x\,\d y\,\d z
\label{eq:0.9}
\end{equation}
or the Gutnik--Nesterenko series \cite{Gu79,Ne96b}
\begin{subequations}
\label{GN}
\begin{equation}
u_n\zeta(3)-v_n=-\frac12\sum_{\nu=1}^\infty\frac{\d}{\d t}
\biggl(\frac{(t-1)(t-2)\dotsb(t-n)}
{t(t+1)(t+2)\dotsb(t+n)}\biggr)^2\bigg|_{t=\nu}.
\label{eq:0.10}
\end{equation}
The fact that each of these representations defines the same $u_n$ and $v_n$ is already a chain of nontrivial analytical
identities. The integral \eqref{eq:0.9} reminds an experienced hypergeometer of the Euler--Pochhammer integral
for generalized hypergeometric functions. The series \eqref{eq:0.10} can be in turn recognised as a hypergeometric
(Barnes-type) integral:
\begin{align}
&
-\frac12\sum_{\nu=1}^\infty\frac{\d}{\d t}
\biggl(\frac{(t-1)(t-2)\dotsb(t-n)}{t(t+1)(t+2)\dotsb(t+n)}\biggr)^2\bigg|_{t=\nu}
\nonumber\\ &\quad
=\frac1{2\pi i}\int_{-i\infty}^{i\infty}\biggl(\frac\pi{\sin\pi t}\biggr)^2
\biggl(\frac{(t-1)(t-2)\dotsb(t-n)}{t(t+1)(t+2)\dotsb(t+n)}\biggr)^2\d t
\nonumber\\ &\quad
=\frac1{2\pi i}\int_{-i\infty}^{i\infty}\frac{\Gamma(t)^4\Gamma(n+1-t)^2}{\Gamma(n+1+t)^2}\,\d t
\label{eq:0.10a}
\end{align}
\end{subequations}
(see \cite{Ne96b} for details). The fact that the right-hand side in \eqref{eq:0.10} indeed represents a linear form in $\zeta(3)$ and $1$
with rational coefficients is quite elementary and uses the partial fraction decomposition \cite{Ne96b,Zu09} of the regular rational function
$$
R(t)=\biggl(\frac{(t-1)(t-2)\dotsb(t-n)}{t(t+1)(t+2)\dotsb(t+n)}\biggr)^2.
$$
This strategy also allows one to gain an explicit arithmetic information about the sequences $\{u_n\}$ and $\{v_n\}$, namely,
$u_n\in\mathbb Z$ (something that can be seen from \eqref{eq:0bin} as well) and $2d_n^3v_n\in\mathbb Z$, where $d_n$ denotes
the least common multiple of the first $n$ positive integers. (The sequence $d_n$ certainly belongs to analytic number theory:
the prime number theorem asserts that $d_n^{1/n}\to e$ as $n\to\infty$.) The same arithmetic can be demonstrated on the basis
of the integrals \eqref{eq:0.9} but based on a different argument\,---\,see \cite{Be79}.

Furthermore, any of representations \eqref{eq:0.9} and \eqref{GN} can be used for estimating the growth of $u_n\zeta(3)-v_n$ as $n\to\infty$.
For example, from \eqref{eq:0.9} we deduce that
\begin{align*}
0<u_n\zeta(3)-v_n
&<\biggl(\max_{[0,1]^3}\frac{x(1-x)y(1-y)z(1-z)}{1-(1-xy)z}\biggr)^n\cdot\frac12\iiint\limits_{[0,1]^3}\frac{\d x\,\d y\,\d z}{1-(1-xy)z}
\\
&=(\sqrt2-1)^{4n}\cdot\zeta(3)
\end{align*}
essentially as an exercise in calculus, while the (more advanced) saddle-point method applied to the single integral in \eqref{eq:0.10a} results in the asymptotics
$$
u_n\zeta(3)-v_n=cn^{-3/2}\,(\sqrt2-1)^{4n}\,(1+O(n^{-1}))
\quad\text{as}\; n\to\infty
$$
for some explicit $c>0$.

Finally, to draw conclusions about the arithmetic of $\zeta(3)$, we assume that it is rational, $p/q$ say, and consider the sequence of then integers
$$
r_n=2d_n^3\,(pu_n-qv_n)=2qd_n^3\,(u_n\zeta(3)-v_n), \quad\text{where}\; n=1,2,\dotsc.
$$
Thus, on one hand $r_n>0$ implying $r_n\ge1$ because of the integrality of the numbers, while on the other hand
$r_n^{1/n}\le e^3(\sqrt2-1)^4(1+o(1))<0.6$ for sufficiently large $n$. The two excluding estimates imply that $\zeta(3)$ cannot be rational,
and the approximations constructed to the number further allow us to measure its irrationality quantitatively:
for any
\begin{equation*}
\mu>\mu_0=1+\frac{4\log(\sqrt2+1)+3}{4\log(\sqrt2+1)-3}
=13.417820\dots,
\end{equation*}
there are only finitely many solutions of the inequality $|\zeta(3)-p/q|<q^{-\mu}$ in integers $p,q$.
(In diophantine approximation theory, if $\alpha$ is an irrational real number then the infimum of $\mu$, for which the inequality
$|\alpha-p/q|<q^{-\mu}$ has at most finitely many solutions, is called the irrationality exponent of $\alpha$ and denoted $\mu(\alpha)$.
Dirichlet's theorem asserts that $\mu(\alpha)\ge2$ for all $\alpha\in\mathbb R\setminus\mathbb Q$; furthermore, considerations
of metric number theory imply $\mu(\alpha)=2$ for almost all real $\alpha$, so that 2 is a typical irrationality exponent.
What is said above is that $\mu(\zeta(3))\le\mu_0$, somewhat insufficiently sharp from a metric point of view, but better than nothing.)

Though we have outlined a proof of the irrationality of $\zeta(3)$ and there are indeed some hypergeometrically looking series and integrals
for constructing Ap\'ery's rational approximations to the number, so far there is no clue to how all this is related to Bailey's transformation~\eqref{7F6tr}.
But we get closer: around 1999, in an unpublished note, Keith Ball gave a different series of rational approximations to $\zeta(3)$, namely,
\begin{subequations}
\label{Ball}
\begin{equation}
\wt u_n\zeta(3)-\wt v_n=n!^2\sum_{\nu=1}^\infty\Bigl(t+\frac n2\Bigr)
\frac{(t-1)\dotsb(t-n)\cdot(t+n+1)\dotsb(t+2n)}{t^4(t+1)^4\dotsb(t+n)^4}\bigg|_{t=\nu}.
\label{eq:0.11}
\end{equation}
He used the clear symmetry $\wt R(-t-n)=-\wt R(t)$ of the rational summand
$$
\wt R(t)=n!^2\Bigl(t+\frac n2\Bigr)\frac{(t-1)\dotsb(t-n)\cdot(t+n+1)\dotsb(t+2n)}{t^4(t+1)^4\dotsb(t+n)^4}
$$
and the related partial-fraction decomposition of $R(t)$ to show that the sum on the right-hand side of \eqref{eq:0.11},
\emph{a priori} living in $\mathbb Q\zeta(4)+\mathbb Q\zeta(3)+\mathbb Q\zeta(2)+\mathbb Q$, has vanishing coefficients of $\zeta(4)$ and $\zeta(2)$.
Ball's proof implied the weaker arithmetic properties $2d_n\wt u_n,2d_n^4\wt v_n\in\mathbb Z$ than those for the representations above
but he could show that $(\wt u_n\zeta(3)-\wt v_n)^{1/n}\to(\sqrt2-1)^4$ as $n\to\infty$ using Stirling's formula for the gamma function
as the hardest ingredient. Unfortunately, because of $e^4(\sqrt2-1)^4>1$, where $e^4$ corresponds to the growth of $d_n^4$,
there seemed to be no way to adapt Ball's construction to a new irrationality proof. Ball did not make explicit the fact that the value of his series \eqref{eq:0.11}
is exactly the same as the one coming from, say, the Gutnik--Nesterenko series~\eqref{eq:0.10}, though he was aware of this on the basis of computation of
a couple of terms and of the asymptotics he obtained. The fact that
$$
u_n\zeta(3)-v_n=\wt u_n\zeta(3)-\wt v_n
$$
was first established in the thesis~\cite{Ri01} of Rivoal, who used (a version of) the Gosper--Zeilberger algorithm of creative telescoping \cite{PWZ97} to verify that
the hypergeometric series
\begin{align}
&
n!^2\sum_{\nu=1}^\infty\Bigl(t+\frac n2\Bigr)
\frac{(t-1)\dotsb(t-n)\cdot(t+n+1)\dotsb(t+2n)}{t^4(t+1)^4\dotsb(t+n)^4}\bigg|_{t=\nu}
\nonumber\\ &\quad
=\frac{(3n+2)!\,n!^7}{(2n+1)!^5}\,{}_7F_6\biggl(\begin{matrix} 3n+2, \, \frac32n+2, \, n+1, \, n+1, \, n+1, \, n+1, \, n+1 \\
\frac32n+1, \, 2n+2, \, 2n+2, \, 2n+2, \, 2n+2, \, 2n+2 \end{matrix} \biggm| 1 \biggr)
\label{eq:0.11a}
\end{align}
\end{subequations}
satisfies Ap\'ery's recursion \eqref{eq:0.3} and the initial data agree with \eqref{eq:0.4}.

The coincidence of the two representations, \eqref{GN} and \eqref{Ball}, is a particular case of the transformation~\eqref{7F6tr}:
it corresponds to the choice of parameters $a=3n+2$ and $b=c=d=e=f=n+1$. Things could have gone quite differently if a hypergeometer
would come in at an early stage and recognise that the Barnes integral \eqref{eq:0.10a} can be transformed into the very-well-poised ${}_7F_6$
followed by a number theorist who would observe the reason for the coefficients of $\zeta(4)$ and $\zeta(2)$ (``the parasites'' as Rivoal called them,
those known to be irrational) to disappear. The phenomenon of this disappearance can be pushed further to construct linear forms, with rational
coefficients, in odd zeta values only using (very-)well-poised hypergeometric functions evaluated at $z=1$ (and $z=-1$). This program was successfully
carried out by Rivoal in \cite{Ri00,Ri01} and in his joint work \cite{BR01} with Ball: for odd $s>1$, their approximating forms were (essentially) given by
\begin{align}
&
2d_n^{s+1}
\cdot n!^{s+1-2r}\sum_{\nu=1}^\infty\Bigl(t+\frac n2\Bigr)
\frac{\prod_{j=1}^{rn}(t-j)\cdot\prod_{j=1}^{rn}(t+n+j)}{\prod_{j=0}^n(t+j)^{s+1}}\bigg|_{t=\nu}
\nonumber\\ &\quad
\in\mathbb Z\zeta(s)+\mathbb Z\zeta(s-2)+\dots+\mathbb Z\zeta(5)+\mathbb Z\zeta(3)+\mathbb Z,
\label{eq:0.19}
\end{align}
where the auxiliary integral parameter $r<s/2$ is of order $r\sim s/\log^2s$ for large $s$. Then the explicit formulae for the forms in \eqref{eq:0.19}
allow one to compute the asymptotic behaviour of them and their coefficients as $n\to\infty$,
and the final step of estimating the number $\delta(s)$ of linearly independent over $\mathbb Q$ among $\zeta(s),\zeta(s-2),\dots,\zeta(3)$ and~1
from below uses a criterion of Nesterenko. The result is $\delta(s)>\frac13\,\log s$ (and $\frac13$ can be replaced by any constant closer to but smaller than
$1/(1+\log2)$ for sufficiently large $s$).

It is quite remarkable that the class of well-poised hypergeometric functions plays such a special role in establishing that infinitely many of the odd zeta values
are irrational. But the transformation \eqref{7F6tr} provides us with slightly more: one can use the full power of \eqref{7F6tr} to produce a sharper
quantitative irrationality of $\zeta(3)$. In Section~7.5 of his book \cite{B35} Bailey discusses the hypergeometric transformation group, of size 1920, that acts
on the six parameters $a,b,c,d,e,f$ of the hypergeometric functions involved. Using this group and the arithmetic ``permutation group'' method
developed by George Rhin and Carlo Viola in \cite{RV96,RV01} one can prove the estimate $\mu(\zeta(3))\le5.513890\dots$ for the irrationality exponent of~$\zeta(3)$.
This is the result originally proved by Rhin and Viola \cite{RV01} in 2001 by applying their novel techniques to the Beukers-type integrals
\begin{align}
I(h,j,k,l,m,q,r,s)
&=\iiint\limits_{[0,1]^3}\frac{x^h(1-x)^ly^k(1-y)^sz^j(1-z)^q}{(1-(1-xy)z)^{q+h-r}}\,\frac{\d x\,\d y\,\d z}{1-(1-xy)z}
\nonumber\\
&\in\mathbb Z\zeta(3)+\mathbb Q
\label{RVI}
\end{align}
that generalise those in~\eqref{eq:0.9}; here the eight positive parameters are subject to the two relations $j+q=l+s$ and $h+m=k+r$
(the latter one in fact defines the parameter $m$ missing in the integral in~\eqref{RVI}). In order to recognise a permutation group acting
on the set of the eight parameters (that, for example, includes the cyclic permutations of the set $(h,j,k,l,m,q,r,s)$), they designed
several birational transformations of the unit cube preserving the measure of integration in~\eqref{RVI} as well as the form of the integrand.
This group can be recognised as Bailey's hypergeometric group from \cite[Section~7.5]{B35} with the help of the identity
\begin{align*}
&
I(b-1,c-1,d-1,a-b-e,a-b-c,
\\[-1.2mm] &\qquad
1+2a-b-c-d-e-f,a-c-d,a-d-f)
\\[1.2mm] &\quad\qquad
=\frac{\splitfrac{\Gamma(1+a)\,\Gamma(b)\,\Gamma(c)\,\Gamma(d)\,\Gamma(1+a-c-d)\,\Gamma(1+a-b-e)}{\times\Gamma(1+a-b-c)\,\Gamma(1+a-d-f)}}
{\Gamma(1+a-b)\,\Gamma(1+a-c)\,\Gamma(1+a-d)\,\Gamma(1+a-e)\,\Gamma(1+a-f)}
\\ &\quad\qquad\quad\times
{}_7F_6\biggl(\begin{alignedat}{7} & a, &\,& 1+\tfrac12a, &\,& \quad\; b, &\,& \qquad c, &\,& \qquad d, &\,& \qquad e, &\,& \qquad f \\[-1.2mm]
 & &\,& \quad \tfrac12a, &\,& {\kern-2.5mm}1+a-b, &\,& 1+a-c, &\,& 1+a-d, &\,& 1+a-e, &\,& 1+a-f \end{alignedat} \biggm| 1 \biggr)
\end{align*}
which follows from either \cite[Theorem~2]{Ne03} or \cite[Theorem~5]{Zu03}. In \cite{Zu04} we show that the Rhin--Viola estimate
for $\zeta(3)$ can be obtained directly on using Bailey's transformation \eqref{7F6tr} and that another hypergeometric transformation
can be applied in a similar fashion to obtain the Rhin--Viola estimate $\mu(\zeta(2))\le5.441242\dots$ for the irrationality exponent
of $\zeta(2)$ (hence, of~$\pi^2$). At the top level of (very-well-poised) hypergeometric hierarchy there is a transformation
of the Barnes-type integral that decomposes into a linear combination of two very-well-poised balanced ${}_9F_8$. This is given
by Bailey in \cite[Section 6.8]{B35} and the corresponding hypergeometric group of order 51840 can be used for estimating
the irrationality exponent of $\zeta(4)$\,---\,the details can be found in~\cite{Zu03} (together with a subtle, and still open, ``denominator
conjecture'' which is required to make the result unconditional). The papers \cite{Zu03,Zu04} are already expository enough to
follow the circle of ideas and arithmetic ingredients around the irrationality of zeta values but one can also check with the reviews
\cite{Fi04,Zu11} for a development of the topic from a broader perspective.

Are there deep reasons for hypergeometric identities and irrationality investigations to be related?
The philosophy is that behind any hypergeometric transformation there is some interesting arithmetic, and one further illustration of the principle
is our recent work \cite{Zu14} in which we prove the record bound $\mu(\zeta(2))\le5.095411\dots$ for the irrationality exponent
of $\zeta(2)$. The transformation used for the proof relates two Barnes-type integrals,
\begin{align}
&
\frac1{2\pi i}\int_{-i\infty}^{i\infty}
\frac{\Gamma(a+t)\,\Gamma(b+t)\,\Gamma(e+t)\,\Gamma(f+t)}{\Gamma(1+t)\,\Gamma(1+a-e+t)\,\Gamma(1+a-f+t)\,\Gamma(g+t)}
\biggl(\frac\pi{\sin\pi t}\biggr)^2\d t
\nonumber\\ &\;
=(-1)^{a+b+e+f}\frac{\Gamma(e+f-a)\,\Gamma(e)\,\Gamma(f)}{\Gamma(g-b)}
\nonumber\\ &\;\quad\times
\frac1{2\pi i}\int_{-i\infty}^{i\infty}
\frac{\Gamma(a-b+g+2t)\,\Gamma(a+t)\,\Gamma(e+t)\,\Gamma(f+t)}{\Gamma(1+a+2t)\,\Gamma(1+a-b+t)\,\Gamma(e+f+t)\,\Gamma(g+t)}
\,\frac\pi{\sin2\pi t}\,\d t,
\label{bmiss}
\end{align}
and is expected to be true for a generic set of \emph{integral} parameters $a,b,e,f,g$; it is only proved in \cite{Zu14}
for a particular (required) set of parameters that depend on a single integral parameter $n$ by application of (a version of)
the Gosper--Zeilberger algorithm of creative telescoping. Both sides of \eqref{bmiss} represent a linear form in $\zeta(2)$ and 1
with rational coefficients. Even more, we expect the companion transformation, in which $(\pi/\sin\pi t)^2$ and $\pi/\sin2\pi t$ are replaced with
$\pi^3\cos\pi t/\allowbreak(\sin\pi t)^3$ and $(\pi/\sin\pi t)^2$, respectively, to be true as well;
the corresponding integrals in that case represent a linear form in~$\zeta(3)$ and 1, with the same coefficient of $\zeta(3)$
as the coefficient of $\zeta(2)$ in the former linear form. The coincidence of the leading coefficients is known to be true in general
thanks to Whipple's transformation \cite[Section~4.5, eq.~(1)]{B35},
\begin{multline}
{}_4F_3\biggl(\begin{matrix} f, \, 1+f-h, \, h-a, \, -N \\
h, \, 1+f+a-h, \, g \end{matrix}\biggm|1\biggr)
=\frac{(g-f)_N}{(g)_N}
\\ \times
{}_5F_4\biggl(\begin{alignedat}{5}
a&, \, & -N&, \, & 1+f-g&, \, & \tfrac12f&, \, & \tfrac12f+\tfrac12 \\
&& h&, \, & 1+f+a-h&, \, & \tfrac12(1+f-N-g)&, \, & \tfrac12(1+f-N-g)+\tfrac12 \\
\end{alignedat}\biggm|1\biggr),
\label{whipple}
\end{multline}
where $N$ is a positive integer, so that the both hypergeometric series in \eqref{whipple} terminate.

Identity \eqref{bmiss} and its companion should be a special case of a hypergeometric-integral identity valid for
generic \emph{complex} parameters. We failed to trace this more general identity in the literature,
though there are a few words about it at the end of Bailey's paper~\cite{B32}:
\begin{quote}
\small
``The formula (1.4)\footnote{This is equation \eqref{whipple} above.}
and its successor are rather more troublesome to generalize, and the final result was unexpected.
The formulae obtained involve five series instead of three or four as previously obtained.
In each case two of the series are nearly-poised and of the second kind, one is nearly-poised
and of the first kind, and the other two are Saalsch\"utzian in type.
In the course of these investigations some integrals of Barnes's type
are evaluated analogous to known sums of hypergeometric series. Considerations of space,
however, prevent these results being given in detail.''
\end{quote}
It is quite similar in spirit to Fermat's famous ``I have discovered a truly marvelous proof of this,
which this margin is too narrow to contain'', isn't it?
Interestingly enough, the last paragraph in Chapter~6 of Bailey's book \cite{B35} again reveals
no details about the troublesome generalization. Did Bailey possess the identity?


\section{Appell's hypergeometric functions\\ and generating functions of Legendre polynomials}
\label{sun}

In the short note \cite{B48} Bailey's gives an elegant generalization of Jacobi's elliptic integral
$$
K(k)=\int_0^{\pi/2}\frac{\d\theta}{\sqrt{1-k^2\sin^2\theta}}, \quad\text{where}\; |k|<1.
$$
Namely, he proves that the two variable extension
$$
I(k,l)=\int_0^{\pi/2}{\kern-3.6mm}\int_0^{\pi/2}\frac{\d\theta\,\d\lambda}{\sqrt{1-k^2\sin^2\theta-l^2\sin^2\lambda}},
\quad\text{where}\; k^2+l^2<1,
$$
can be evaluated in terms of the elliptic integrals as follows:
$$
I(k,l)=\frac2{1+l'}\,K(k_1)\,K(k_2),
$$
where
$$
k_1=\frac{k'-\sqrt{1-k^2-l^2}}{1+l'}, \quad
k_2=\frac{\sqrt{(1+k)(l'+k)}-\sqrt{(1-k)(l'-k)}}{1+l'}
$$
and $k',l'$ are the moduli complementary to $k,l$ (in other words, $k'=\sqrt{1-k^2}$ and $l'=\sqrt{1-l^2}$).

The beauty of the evaluation is mainly in the ingredients of the proof which uses the transformation~\cite{AKdF26}
\begin{align}
&
F_2\biggl(a;\,a-b+\tfrac12,b;\,c,2b\biggm|\frac{X}{(1+Y)^2},\frac{4Y}{(1+Y)^2}\biggr)
\nonumber\\ &\quad
=(1+Y)^{2a}F_4(a,a-b+\tfrac12;\,c,b+\tfrac12\mid X,Y^2)
\label{App}
\end{align}
of Appell's hypergeometric functions in two variables
\begin{align*}
F_2(a;\,b_1,b_2;\,c_1,c_2\mid x,y)
&=\sum_{m,n\ge0}\frac{(a)_{m+n}(b_1)_m(b_2)_n}{m!\,n!\,(c_1)_m(c_2)_n}\,x^my^n
\\ \intertext{and}
F_4(a,b;\,c_1,c_2\mid x,y)
&=\sum_{m,n\ge0}\frac{(a)_{m+n}(b)_{m+n}}{m!\,n!\,(c_1)_m(c_2)_n}\,x^my^n,
\end{align*}
an Euler-type double integral for $F_2$ and a special reduction of $F_4$ to single-variable hypergeometric functions.
And the latter special reduction requires our special attention, because it generated several applications in different parts
of analysis and was a crucial part of proofs of certain number-theoretical identities. It is
\begin{equation}
F_4(a,b;\,c,a+b-c+1\mid X(1-Y),Y(1-X))
={}_2F_1\biggl(\begin{matrix} a, \, b \\ c \end{matrix}\biggm|X\biggr)
\,{}_2F_1\biggl(\begin{matrix} a, \, b \\ a+b-c+1 \end{matrix}\biggm|Y\biggr),
\label{F4}
\end{equation}
valid inside simply-connected regions surrounding $X=0$, $Y=0$ for which
$$
|X(1-Y)|^{1/2}+|X(1-Y)|^{1/2}<1.
$$
This reduction is quite different from those given earlier in \cite{AKdF26} as the result is a product
of two Euler--Gauss hypergeometric functions rather than a single one of the type. Formula \eqref{F4} was published
in \cite{B33,B34} and shortly thereafter in the book \cite{B35} and it is clear from the comments in the latter
as well as in its use in later works of Bailey that this was one of his personal favourites.

This time the number theory counterpart came in 2011 from China, where motivated by Ramanujan's beautiful formulae and systematically
experimenting with binomial expressions Z.-W.~Sun observed \cite{Su11} some that produce approximations to simple multiples of $1/\pi$.
Some typical examples of Sun's production were the identities
\begin{equation}
\begin{aligned}
\sum_{n=0}^\infty\frac{7+30n}{(-1024)^n}\,{\binom{2n}{n}}^2T_n(34,1)&=\frac{12}{\pi},
\\
\sum_{n=0}^\infty\frac{2+15n}{972^n}\,\binom{2n}{n}\binom{3n}{2n}T_n(18,6)&=\frac{45\sqrt3}{4\pi},
\end{aligned}
\label{sun-ex}
\end{equation}
where $T_n(b,c)$ denotes the coefficient of $x^n$ in the expansion of $(x^2+bx+c)^n$, explicitly
\begin{equation*}
T_n(b,c)=\sum_{k=0}^{\lfloor n/2\rfloor}\binom{n}{2k}\binom{2k}{k}b^{n-2k}c^k.
\end{equation*}
Although the formulae are not very practical for computing $1/\pi$ (hence, $\pi$ itself), they very much resemble
the formulae that were given by Srinivasa Ramanujan \cite{Ra14} in~1914 such as
\begin{equation*}
\sum_{n=0}^\infty\frac{1103+26390n}{396^{4n}}\,{\binom{2n}{n}}^2\binom{4n}{2n}
=\frac{99^2}{2\pi\sqrt2}
\end{equation*}
which converges to the multiple of $1/\pi$ on the right-hand side very rapidly. In both situations we have the pattern
$$
\sum_{n=0}^\infty(A+Bn)\,u(n)\,z_0^n=\frac C\pi,
$$
where $A$ and $B$ are certain integers, $C$ is an algebraic number and $z_0$ is a rational close to the origin.
In Ramanujan's cases from \cite{Ra14} the series $\sum_{n=0}^\infty u(n)z^n$ is a hypergeometric ${}_3F_2$ series (with a special choice
of the parameters), which clearly does not happen for the examples in~\eqref{sun-ex}.
In fact, the main feature of Ramanujan's identities for $1/\pi$ and their later generalizations is that
the function $\sum_{n=0}^\infty u(n)z^n$ satisfies a third order (arithmetic) linear differential equation
with regular singularities which happens to be a symmetric square of a second order differential equation.
The interested reader is advised to consult the surveys \cite{BBC09,Zu08} on the topic of Ramanujan-type formulae for $1/\pi$
and here we will only indicate a classical hypergeometric instance of such hypergeometric functions $\sum_{n=0}^\infty u(n)z^n$,
which is known as Clausen's identity:
\begin{equation}
{}_3F_2\biggl(\begin{matrix} \frac12, \, r , \, 1-r \\ 1, \, 1 \end{matrix}\biggm|4x(1-x)\biggr)
={{}_2F_1\biggl(\begin{matrix} r, \, 1-r \\ 1 \end{matrix}\biggm|x\biggr)}^2.
\label{clauclau}
\end{equation}
The \emph{arithmetic} cases correspond to the choice $r\in\{\frac12,\frac13,\frac14,\frac16\}$.

The functions $\sum_{n=0}^\infty u(n)z^n$ in Sun's examples from \cite{Su11} satisfy fourth order linear differential equations,
so that the structure of those identities is somewhat different from the one in Ramanujan's situations.
But they can be cast in a more hypergeometric form, because the binomial sums $T_n(b,c)$ are recognised as
$$
T_n(b,c)=(b^2-4c)^{n/2}P_n\biggl(\frac b{(b^2-4c)^{1/2}}\biggr),
$$
where
$$
P_n(x)={}_2F_1\biggl(\begin{matrix} -n, \, n+1 \\ 1 \end{matrix}\biggm| \frac{1-x}2 \biggr)
$$
are the classical Legendre polynomials. Then the formulae on Sun's list read
$$
\sum_{n=0}^\infty(A+Bn)\,\frac{(r)_n(1-r)_n}{n!^2}\,P_n(x_0)\,z_0^n=\frac C\pi
$$
with some algebraic $x_0,z_0$ and, as before, $r\in\{\frac12,\frac13,\frac14,\frac16\}$. The series
rang a bell with us and after a little search in the literature about generating functions of orthogonal polynomials we found that
\begin{equation}
\sum_{n=0}^\infty\frac{(r)_n(1-r)_n}{n!^2}P_n(x)z^n
={}_2F_1\biggl(\begin{matrix} r, \, 1-r \\ 1 \end{matrix}\biggm| \frac{1-\rho-z}2 \biggr)
\cdot{}_2F_1\biggl(\begin{matrix} r, \, 1-r \\ 1 \end{matrix}\biggm| \frac{1-\rho+z}2 \biggr),
\label{e06}
\end{equation}
where $\rho=\rho(x,z):=(1-2xz+z^2)^{1/2}$, as a particular entry on the rich list \cite{Br51} of findings
of Fred Brafman in his 1951 thesis. One can notice that the two identities \eqref{clauclau} and \eqref{e06} have a very
similar shape and the only difference is that the product of two different specializations of ${}_2F_1$ in \eqref{e06}
is replaced with the square of a single specialization. It has become a quite enjoyable adventure for three of us,
Heng Huat Chan, James Wan and myself, to adapt in \cite{CWZ13} the methods of proofs of Ramanujan\break\hbox{(-type)} formulae for $1/\pi$
to the new settings and rigorously establish the experimental observations of Sun from~\cite{Su11}.

The fact that Brafman's generating function \eqref{e06} looks very much like Bailey's formula \eqref{F4} is not surprising
as Brafman used the latter one for derivation of the former by using
$$
F_4(r,1-r;\,1,1\mid X(1-Y),Y(1-X))
=\sum_{n=0}^\infty\frac{(r)_n(1-r)_n}{n!^2}P_n\biggl(\frac{X+Y-2XY}{Y-X}\biggr)(Y-X)^n.
$$
In fact, Brafman's identity in \cite{Br51} was a general generating functions of Jacobi polynomials, something that Bailey
had found himself earlier in \cite{B38} (but the trigonometric way Bailey stated his result obscured its applicability).
Brafman went much further in his dedication to generating functions of orthogonal polynomials, and in our joint project \cite{WZ12}
with Wan we used another (hypergeometric!) theorem of Brafman from \cite{Br59} and its generalization given by H.\,M.~Srivastava \cite{Sr75}
to produce some new generating functions of rarefied Legendre polynomials, valid in a neighbourhood of $X=Y=1$:
\begin{align*}
&
\sum_{n=0}^\infty \frac{(\frac12)_n^2}{n!^2} \,P_{2n}\biggl(\frac{(X+Y)(1-XY)}{(X-Y)(1+XY)}\biggr)
\biggl(\frac{X-Y}{1+XY}\biggr)^{2n}
\\ &\qquad
=\frac{1+XY}{2}
\,{}_2F_1\biggl(\begin{matrix} \frac12, \, \frac12 \\ 1 \end{matrix}\biggm| 1-X^2 \biggr)
\,{}_2F_1\biggl(\begin{matrix} \frac12, \, \frac12 \\ 1 \end{matrix}\biggm| 1-Y^2 \biggr)
\end{align*}
and
\begin{align*}
&
\sum_{n=0}^\infty \frac{(\frac13)_n(\frac23)_n}{n!^2}
\,P_{3n}\biggl(\frac{X+Y-2X^2Y^2}{(X-Y)\sqrt{1+4XY(X+Y)}}\biggr)
\biggl(\frac{X-Y}{\sqrt{1+4XY(X+Y)}}\biggr)^{3n}
\\ &\qquad
=\frac{\sqrt{1+4XY(X+Y)}}3
\,{}_2F_1\biggl(\begin{matrix} \frac13, \, \frac23 \\ 1 \end{matrix}\biggm| 1-X^3 \biggr)
\,{}_2F_1\biggl(\begin{matrix} \frac13, \, \frac23 \\ 1 \end{matrix}\biggm| 1-Y^3 \biggr),
\end{align*}
and use them for proving other observations of Sun~\cite{Su11}.
In showing the equalities we needed to establish a certain generalization of Bailey's identity \eqref{F4}
to non-hypergeometric settings, when the factors
$$
{}_2F_1\biggl(\begin{matrix} a, \, b \\ c \end{matrix}\biggm|z\biggr)
$$
in the product \eqref{F4} are replaced with solutions $\sum_{n=0}^\infty u(n)z^n$ of more general (arithmetic)
differential equation of order 2 with four regular singularities. We refer the interested reader to our original paper~\cite{WZ12}
for the details of the formula and story.

\medskip
There is a more significant difference between Clausen's formula \eqref{clauclau} and Bailey's \eqref{F4}:
the former one depends on a single variable while there are two variables involved in the latter.
There are efficient algorithms to determine whether a given solution $F(z)$ of a third order (Picard--Fuchs) linear differential equation
can be written in the form $F(z)=\alpha(z)\cdot f(\beta(z))^2$ where $\alpha(z)$ and $\beta(z)$ are algebraic functions
and $f(z)$ satisfies a second order linear differential equation. At the same time, no general algorithm is known
to write for a given function $F(X,Y)$  a representation
\begin{equation}
F(X,Y)=\alpha(X,Y)\cdot f(\beta(X,Y))f(\gamma(X,Y))
\label{e0}
\end{equation}
 (or even to check whether one exists) with some \emph{algebraic} functions $\alpha(X,Y)$, $\beta(X,Y)$ and $\gamma(X,Y)$ and $f(z)$ a solution
of a second order (Picard--Fuchs) differential equation. Many cases when we expect such factorisations \eqref{e0} to exist
are suggested by later experimental findings of Sun, which he fed in the later updates of~\cite{Su11}. For example,
we expect factorisations of the form~\eqref{e0} to hold for the generating functions
\begin{equation*}
\sum_{n=0}^\infty P_n(x)^3z^n \quad\text{and}\quad
\sum_{n=0}^\infty\binom{2n}nP_n(x)^2z^n,
\end{equation*}
where $P_n(x)$ are again the Legendre polynomials, but we only know \cite{Zu14b} a formula for the latter one when $z=4(x^2-1)/(x^2+3)^2$.

Most of the known factorisations \eqref{e0} follow from Bailey's formula \eqref{F4} or from its generalization in \cite{WZ12}.
Aware of our interest in collecting such examples, Frits Beukers communicated to us in 2013 his personal finding of factorisation of Appell's
function $F_2$, which he had come across accidentally when studying reducible cases of GKZ hypergeometric functions and proved
using somewhat \emph{ad hoc} methods. Beukers' identity is
\begin{align*}
&
F_2\biggl(a+b-\tfrac12;\,a,b;\,2a,2b\biggm|\frac{4u(1-u)(1-2v)}{(1-2uv)^2},\frac{4v(1-v)(1-2u)}{(1-2uv)^2}\biggr)
\\ &\quad
=(1-2uv)^{-1+2a+2b}
{}_2F_1\biggl(\begin{matrix} a+b-\tfrac12, \, a \\ 2a \end{matrix}\bigg|4u(1-u)\biggr)
\,{}_2F_1\biggl(\begin{matrix} a+b-\tfrac12, \, b \\ 2b \end{matrix}\bigg|4v(1-v)\biggr),
\end{align*}
and it indeed follows from the combination of \eqref{App} with~\eqref{F4}
(Bailey himself wrote an equivalent form of Beukers' formula in his 1938 paper~\cite{B38}; see eq.~(3.1) there).

\medskip
Bailey's last publication \cite{B59} in 1959 was again about the formula \eqref{F4}, this time from a historical perspective.
While working on the obituary notice \cite{B54} of Ernest William Barnes for the London Mathematical Society, Bailey was given
access to some unpublished materials from Mrs.\ Barnes. One manuscript prepared by Barnes in 1907 for publication (but never published)
gave the formula and its proof. Bailey writes in \cite{B59}:
\begin{quote}
\small
``The formula (1.1)\footnote{This is equation \eqref{F4} above.} was obtained by myself in 1933, and I did not know until the discovery of the manuscript that Barnes had
obtained the formula a quarter of a century before I did. The manuscript was quite a lengthy one and gave different forms of (1.1) with
a number of diagrams illustrating the regions of validity in the different cases. I sent the work to Professor Watson, and at one time he contemplated
writing a short note on it, but was deterred by the fact that the state of printing at the time was very difficult.''
\end{quote}

We can speculate that such finding was a personal tragedy to Bailey and one possible reason for the decline of his research activities.
It is certain to us that the formula \eqref{F4} is Bailey's, and we admire his act of honesty in reporting on the discovered manuscript of Barnes in public.

\medskip
We conclude this section by presenting one more episode on the function $F_4$.
In his brief exposition about Appell's hypergeometric functions \cite[Chapter~9]{B35}, Bailey notices (in Section~9.3, after
recording Euler-type double integrals for the functions $F_1$, $F_2$ and $F_3$) that ``[t]here appears
to be no simple integral representation of this type for the function $F_4$.'' A formula of this type was given by Burchnall and Chaundy in 1940
in the paper~\cite{BC40} where they introduced a symbolic notation for expressing Appell's functions by means of
single-variable hypergeometric functions. Bailey discusses the results from \cite{BC40} and gives his own proof of
the double integral representation of $F_4$ in~\cite{B41a}. Surprisingly enough, this discovery is not reflected
in later Slater's book~\cite{Sl66}, where she misleadingly indicates (in Section~8.2) ``no similar integral for $F_4$ has been found.''

\section{An algebraic identity: Variations on $q$-analogues of zeta values}
\label{qsect}

An \emph{algebraic}, as appears in the title of Bailey's short communication \cite{B36c}, identity in question
is in fact of combinatorial nature and the adjective `algebraic' is probably used as a synonym to `non-analytical'
(though the identity is analytical as well):
\begin{equation}
\sum_{n=1}^\infty\frac{q^n}{(1-q^n)^2}\sum_{l=1}^n\frac1{1-q^l}
=\sum_{n=1}^\infty\frac{n^2q^n}{1-q^n}.
\label{id-Bell}
\end{equation}
He learned about it from Hardy and attributed it to E.\,T.~Bell;
in \cite{B36c} Bailey records three different proofs of the identity.
If we denote $\sigma_k(n)=\sum_{d\mid n}d^k$
the sum of the $k$th powers of the divisors of~$n$ then the Lambert series on the right-hand side
of \eqref{id-Bell} can be alternatively written as $\sum_{n=1}^\infty\sigma_2(n)q^n$.
As Bell mentions himself \cite{Be35} (see the related discussion in the later note \cite{B36d} of Bailey)
the identity has the following combinatorial interpretation
stated without proof by Liouville~\cite{Li67}:
the number of representation of a positive integer $n$ in the form $ab+bc+cd+de$, where $a,b,d,e>0$
and $c\ge0$ are integers, is $\sigma_2(n)-n\sigma_0(n)$.

As pointed out in the introduction, the identity is a $q$-analogue of Euler's $\zeta(2,1)=\zeta(3)$ and it was stated by us
with a mistake in the survey \cite{Zu03a} on the so-called multiple zeta values (MZVs). (At the time of writing \cite{Zu03a} we were
not aware of Bailey's \cite{B36c} and its predecessors.) The mistake was soon after corrected by D.~Bradley in \cite{Br05}
and the whole area of $q$-MZVs has exploded in the last years; we limit ourselves here to mention of the two excellent representatives
of the explosion\,---\,Bachmann's bi-brackets and multiple Eisenstein series in \cite{Ba15}
and the structural relations of $q$-MZVs in \cite{CEM15} by J.~Castillo Medina, K.~Ebrahimi-Fard and D.~Manchon.

The result \eqref{id-Bell} is particularly interesting as it gives two representations of a generating
function linked with the world of Maass forms~\cite{LZ01} rather than modular forms. The latter fact makes
this story quite disjoint with that for the Rogers--Ramanujan identities \cite{Ro17,RR19} when the both sides
of such an identity are expected to represent some modular forms. (Recent experimental discoveries
of Kanade and Russell in \cite{KR15} give some evidence in that the modularity may not be always a feature.)
This makes Slater's remark \cite{Sl62} about connection of the two related equations from \cite{B36c},
\begin{align}
(1-q)\sum_{n=1}^\infty\frac{(1-q^{2n+1})q^n}{(1-q^n)^2(1-q^{n+1})^2}
\biggl(\frac{1+q}{1-q}+\frac{1+q^2}{1-q^2}+\dots+\frac{1+q^n}{1-q^n}\biggr)
&=\sum_{n=1}^\infty\frac{n^2q^n}{1-q^n}
\label{id-Bell2}
\\ \intertext{and}
(1+q)\sum_{n=1}^\infty\frac{(1+q^{2n+1})q^n}{(1-q^n)^2(1+q^{n+1})^2}
\biggl(\frac{1+q^2}{1-q^2}+\frac{1+q^4}{1-q^4}+\dots+\frac{1+q^{2n}}{1-q^{2n}}\biggr)
&=\sum_{n=1}^\infty\frac{n^2q^n}{1-q^n},
\label{id-Bell3}
\end{align}
to the  Rogers--Ramanujan identities unjustified, though indeed the left-hand sides of both \eqref{id-Bell2} and \eqref{id-Bell3}
smell like the sum-parts of some Rogers--Ramanujan-type identities. Apart from this similarity and the similarity of the methods used
in proofs, they share little.

We can also remark that the non-modularity of the series
$$
\zeta_q(s)=\sum_{n=1}^\infty\frac{n^{s-1}q^n}{1-q^n}
$$
for odd integers $s\ge1$ and the fact that $(1-q)^s\zeta_q(s)\to(s-1)!\,\zeta(s)$ for any integer $s>1$ make the series perfect $q$-analogues
of (odd) $q$-zeta values \cite{Zu02b}. The first arithmetic result in their direction was obtained by Paul Erd\H os \cite{Er48} in 1948,
who proved the irrationality of the $q$-harmonic series $\zeta_q(1)$ for $q$ a reciprocal of an integer $>1$.
The analogies between the $q$-series $\zeta_q(s)$ and quantities $\zeta(s)$ for positive $s>1$ have deep hypergeometric roots;
in particular, the $q$-basic version of the very-well-poised construction from Section~\ref{zeta} produces the linear forms in odd $q$-zeta values
$\zeta_q(s)$ with coefficients from $\mathbb Q(q)$. This was exploited in our joint paper \cite{KRZ06} with Krattenthaler and Rivoal to estimate the dimension
of the $\mathbb Q$-space spanned by the odd $q$-zeta values under the conditions on $q$ similar to that of Erd\H os in~\cite{Er48}.
This is a number-theoretical $q$-analogue of Rivoal's theorem \cite{BR01,Ri00,Ri01}.

\section{The Erd\H os number of W.\,N.~Bailey}
\label{erdos}

Bailey's lack of collaboration is striking:
among the 75 papers and one book authored by him (which are carefully listed by Slater in~\cite{Sl62}) there is only one work coauthored.
This is the joint paper \cite{WB38} with John Macnaghten Whittaker published in 1938, which is one page long and places Bailey second on the authors' list!
Already this single little publication makes Bailey's Erd\H os number (that is, his collaboration distance to Erd\H os) finite, namely equal to~4, in view of the collaboration chain
J.\,M.~Whittaker--R.~Wilson (\emph{J.\ London Math.\ Soc.}\ \textbf{14} (1939), 202--208),
R.~Wilson--A.\,J.~Macintyre	(5 joint papers recorded by the \emph{MathSciNet} including, e.g., \emph{Proc.\ London Math.\ Soc.}\ (2) \textbf{47} (1940), 60--80)
and A.\,J.~Macintyre--P.~Erd\H os (\emph{Proc.\ Edinburgh Math.\ Soc.}\ (2) \textbf{10} (1954), 62--70).
The sequence was communicated to us by Jonathan Sondow. Before his sending we believed that Bailey's Erd\H os number
was infinite, and designed a different path from Bailey to Erd\H os through J.\,M.~Whittaker's biological father, Edmund Taylor Whittaker (who clearly had had some mathematical impact on his son).
According to the \emph{MathSciNet} (accessed on 30 October 2016), the Erd\H os number of E.\,T.~Whittaker is equal to 4 and it goes through his famous
collaboration \emph{A course of modern analysis} with George Neville Watson (with the latter placed second!),
who coauthored the memoir \cite{ABBW85} with Bruce Berndt, who in turn had a joint paper with Sarvadaman Chowla, an Erd\H os's collaborator.
(The search in \emph{Zentralblatt} reveals a different collaboration path from Whittaker the father to Erd\H os, through
Harry Bateman, Stephen Rice and Nicolaas de Bruijn\,---\,not a bad one either!)
The fact that Watson authored a paper some 20 years after passing away immediately caused a question of whether another G.\,N.~Watson exists,
with a similar circle of mathematical interests. Here we reproduce the response of Bruce Berndt (dated 30 October 2016) to our query:
\begin{quote}
\small
``In regard to our AMS Memoir on Chapter 16 of the second notebook,
coauthored with Watson, the three living authors felt that Watson should
be given the credit that he deserved, and so we listed him as a coauthor.
Fortunately, the Memoir editors were sympathetic with our view.  I tried
at least one other time to list Watson as a coauthor, namely for my paper
with Ron Evans, Extensions of asymptotic expansions from Chapter~15 of
Ramanujan's second notebook, \emph{J.~Reine Angew.~Math.} \textbf{361} (1985),
118--134.  However, the editors of \emph{Crelle} told us that each author needed
to give explicit agreement for publication.  Since Watson was dead,
it was of course impossible to get his consent.''
\end{quote}
We can only add to the story that these days consents for publication are even tougher.

\medskip
One remarkable story of what came out from Bailey's interest in basic hypergeometric functions and generalizations
of the famous Rogers--Ramanujan identities was his correspondence with Freeman Dyson in the 1940s.
Dyson records this in \cite{Dy88} as one of his visits to Ramanujan's garden. We are fortunate to include
copies of those few (originally hand-written) letters from Bailey to Dyson below as an appendix to the paper,
in which one can get a feeling for the personality of W.\,N.~Bailey.
There it becomes transparent that the correspondence has originated a method that would be later published by Bailey in \cite{B47a,B49} and
much later coined the name ``Bailey's lemma'' by George Andrews~\cite{An84} in~1984.
This method and its generalizations have had a great impact on many developments in hypergeometric functions,
combinatorics and number theory~\cite{Wa01}.

\medskip
A somewhat different implication of Bailey's work was on classifying algebraic hypergeometric functions, and our exposition of this here is somewhat approximate
(as we could not witness the events ourselves).
In 1957 Antonius Levelt started his PhD work under supervision of de Bruijn. An initial task for him was to find an appropriate way of representing
the material in Bailey's book \cite{B35}. One should take into account that the late 1950s and 1960s were an explosion of Bourbaki's ideas and abstract algebraic approaches
in mathematics\,---\,the one-page exposition in \cite[Section~8.4]{Re88} is an excellent reference for gaining the ``spirit of the epoch''\,---\,and
Levelt followed Grothendieck's lectures in Paris for two years. The subject of hypergeometric functions was hardly considered as ``sexy'' at that time
and he found difficulties in just appreciating the magical formulae that were always established by a combination of certain tricks. What Levelt managed to do,
however, was finding a different route of his own to the subject and relating the monodromy of hypergeometric differential equations to a simple and elegant
problem in linear algebra. He defended his thesis \cite{Le61} \emph{cum laude} in~1961, and in the last chapter of the thesis he illustrated the power of the
machinery he developed by proving certain particular entries from Bailey's book. Levelt's work was picked up later by Frits Beukers and Gert Heckman \cite{BH89}
to give a complete description of groups that show up as the (Zariski closures of) hypergeometric monodromy groups and, by this means, to give a complete
list of hypergeometric functions~\eqref{eq01} which are algebraic.


\setcounter{footnote}{0}
\renewcommand{\thefootnote}{\fnsymbol{footnote}}
\numberwithin{equation}{section}
\makeatletter
\tagsleft@true
\makeatother

\section*{Appendix: W.\,N.~Bailey's letters to F.~Dyson}

A historical context of the correspondence can be found in the visit entitled ``III.~Bailey'' of Dyson's ``walk through Ramanujan's garden'' \cite{Dy88}.

\hbox to\hsize{\kern6mm\hrulefill\kern6mm}

\hbox to\hsize{\hss\hbox{\vtop{\hsize=36mm
\leftline{17 Prince's Road,}
\leftline{Heaton Moor,}
\leftline{Stockport,}
\leftline{Cheshire.}
\smallskip
\leftline{Dec.~22/43}
}}}

\vskip4mm
\noindent
Dear Dyson, \\
\hbox to 24mm{\hfill}%
I am writing to you in connection with my paper for
the L.M.S. I had a letter from Hardy two days ago, \& he
enclosed your report (a rather unusual thing to do, but
very useful in the circumstances).

Many thanks for your comments. Actually I have not seen
Rogers' papers for some years \& I am in the unfortunate
position of not being able to look them up. The only
copy in Manchester (as far as I am aware) was destroyed
in the blitz in 1940, or I should certainly have refreshed my
memory. I am now trying to purchase the parts of the L.M.S.\
Proceedings in which they appear, but don't know whether I
shall be successful. I shall, of course, after the last part of
\S\,1 and the first part of \S\,2 to meet your criticism.

With regard to the formulae (4.3) \& (5.3), I did not remember
that they were given by Rogers, and I couldn't look up this
paper either. I had, however, discovered them in a paper by
Jackson, a copy of which I enclose as it might interest you.
[He sent me two copies]. Jackson makes no reference
to Rogers. The formulae are the last on p.~175 \& the first
on p.~176 (with a misprint). My (6.3) was, I thought, new, but
I think the third formula on p.~170 of Jackson's paper is meant to
be the same, but it is wrong. I am not sure that he deserves
to be quoted when he simply states the formula \& gives it
incorrectly, but I shall probably put in a reference. Jackson is
terribly careless \& has caused me a good deal of trouble, but
he has a good many curious results in this paper.

With regard to your formulae for products in which the powers of $x$
advance by 27, I am afraid I don't see how you obtained
them. I should very much like to know. The products in the
first three are the products occurring in my formulae for 9's
with $x$ replaced by $x^3$, as I suppose you noticed. If you like, you
could make a short paper about them \& ask for it to follow mine
probably in the Proceedings or write one for the journal quoting
what is necessary of my paper. If you still don't think it worth
while making a separate paper, I will incorporate the formulae
in my paper with due acknowledgements. I should
however have to give at least some indication of how they
were obtained. I may be dense, but they don't seem at all
obvious to me. I shall certainly be interested to know
how you got them.

Yours sincerely,

W.~N.~Bailey

\hbox to\hsize{\kern6mm\hrulefill\kern6mm}
\vskip6mm


\hbox to\hsize{\hss\hbox{\vtop{\hsize=36mm
\leftline{17 Prince's Road,}
\leftline{Heaton Moor,}
\leftline{Stockport.}
\smallskip
\leftline{Dec.~24/43}
}}}

\vskip4mm
\noindent
Dear Dyson, \\
\hbox to 24mm{\hfill}%
Just a line to tell you not to bother writing out proofs
of your identities\,---\,if I am not already too late. I have
rather belatedly found out how you got them.

Two of them come from
\begin{align*}
&
1+\sum_{n=1}^\infty\frac{(-1)^n[ax^3]_{n-1}(1-ax^{6n})a^{4n}x^{\frac12(27n^2-3n)}}{[x^3]_n}
\\ &\qquad\qquad
=\prod_1^\infty(1-ax^n)
\sum_{n=0}^\infty\frac{a^nx^{n^2}[ax^3]_{n-1}}{x_n!\,(ax)_{2n-1}}
\end{align*}
where $[\;\;]$ denote the powers of $x$ advance by 3. Two of your
identities come by taking $a=1$ \& $a=x^3$. I suppose you got
your results in the same way as I have done, by taking
$\alpha_r=0$ unless $r$ is a multiple of 3. I should also have got the
last formula of the original three, but something has gone
wrong \& I cannot find the error at present. Perhaps it is too
near Xmas!

I think several other formulae should be obtainable by similar
methods, but of course they may not be new.

By the way, the third formula on p.~170 of Jackson's paper is
evidently meant to be one of mine. All that is needed to
put it right is to change the sign of $q$ on one side.

If you decide to let your identities go in my paper I will send the
paper on to you when finished. They would put a finishing touch on
my paper, \& definitely increase its value, but I think you would be
perfectly justified in making a separate paper.
Best seasonal greetings.

W.~N.~Bailey

\hbox to\hsize{\kern6mm\hrulefill\kern6mm}
\vskip6mm


\hbox to\hsize{\hss\hbox{\vtop{\hsize=36mm
\leftline{17 Prince's Road,}
\leftline{Heaton Moor,}
\leftline{Stockport,}
\leftline{Cheshire.}
\smallskip
\leftline{Jan.~5/44}
}}}

\vskip4mm
\noindent
Dear Dyson, \\
\hbox to 24mm{\hfill}%
I am writing to give you information this time\,---\,not to worry you.

Your method of getting ``27'' identities appears to be equivalent
to the method I found. I have added two paragraphs to my paper,
the first bringing in the idea of making various $\alpha_r=0$ (or
taking various $\cos r\theta=0$), \& getting in particular two of your ``27''
identities. In the last paragraph I give a list of all the identities
you sent, but I don't attempt to provide proofs. Actually I
have only looked at Rogers' methods from one aspect, \& I
doubt whether proofs of all your identities could be got by the
methods of my paper.

I have now got a copy of Rogers' paper which he wrote in
1917\,---\,borrowed it from Jackson who very kindly sent me some
other including Rogers' ``Third Memoir'', but unfortunately not
the other two.

By the way, the formula for $\displaystyle\prod_1^\infty\biggl(\frac{1-x^{9n}}{1-x^n}\biggr)$
is a brute to get. It is easy enough to get $\underline{\text{a}}$ series for it,
but I find it an awful business to get your series.

I thought it wasn't necessary to worry you again with my
paper, particularly as you are so busy, but I wanted to let
you know that I have put your identities on record.

Yours sincerely,

$\underline{\text{W.~N.~Bailey}}$

\hbox to\hsize{\kern6mm\hrulefill\kern6mm}
\vskip6mm

\hbox to\hsize{\hss\hbox{\vtop{\hsize=36mm
\leftline{17 Prince's Road,}
\leftline{Heaton Moor,}
\leftline{Stockport,}
\leftline{Ches.}
\smallskip
\leftline{Feb.~13/44}
}}}

\vskip4mm
\noindent
Dear Dyson, \\
\hbox to 24mm{\hfill}%
I was interested in your last letter which I received some time ago.
I should think your proof of the identity for $\prod\limits_1^\infty(1-q^{9n})$
is as short as can be expected. I have got $a$-generalisation for it,
but unfortunately it is not at all elegant. It is, in fact,
\begin{multline*}
1+\sum_1^\infty\frac{(-1)^n\{ax^3\}_{n-1}}{{x^3}_n!}
(1-ax^{6n})(1-x^{3n}+ax^{6n})a^{4n-1}x^{\frac12(27n^2-9n)}
\\
=\prod_{m=1}^\infty(1-ax^m)
\sum_{n=0}^\infty\frac{\{ax^3\}_n}{x_n!\,(ax)_{2n+1}}\,a^nx^{n^2+n},
\end{multline*}
where $\{a\}_n=(1-a)(1-ax^3)\dotsb(1-ax^{3n-3})$, \\
and I don't see any nicer way of writing it.

I am really writing to you to let you know how things are going. I have been
pretty busy in other ways lately, but I have done enough to feel rather
disappointed with this sort of thing. After studying Rogers' papers, I was led
to put things in this way:

If $\beta_n=\sum\limits_{r=0}^n\alpha_ru_{n-r}v_{n+r}$, and
$\gamma_n=\sum\limits_{r=n}^\infty\delta_ru_{r-n}v_{r+n}$, then
$\sum\limits_{n=0}^\infty\alpha_n\gamma_n=\sum\limits_{n=0}^\infty\beta_n\delta_n$,
provided of course that convergence conditions are satisfied.
This leads, in particular, to all the known transformations of ordinary
hypergeometric series. In fact, it is substantially equivalent
to the method used in my tract, though I think this is rather a more illuminating
way of putting it. Similarly it is substantially equivalent to the method I used in my
last paper to find the transformation of a nearly poised basic series. One form we can
take is
\begin{align*}
\beta_n&=\sum_{r=0}^n\frac{\alpha_r}{(q)_{n-r}(aq)_{n+r}}
\\
\gamma_n&=\sum_{r=n}^\infty\frac{\delta_r}{(q)_{r-n}(aq)_{r+n}}
\end{align*}
and then
\begin{alignat*}{2}
(aq)_{2n}\beta_n
&=\text{Rogers'}\; m_{2n} &\quad\text{if}\; a&=1
\\
&=\quad\;\;{}''\quad\;\; m_{2n+1} &\quad\text{if}\; a&=q
\end{alignat*}
while
\begin{alignat*}{2}
\delta_n(aq)_{2n}
&=a_{2n} &\quad\text{if}\; a&=1
\\
&=a_{2n+1} &\quad\text{if}\; a&=q
\end{alignat*}
and $\gamma_n=b_{2n}$ or $b_{2n+1}$, $\alpha_n=m_{2n}$ or $m_{2n+1}$.

The formula for $\gamma_n$ gives, if $\delta_r=(\rho_1)_r(\rho_2)_r(aq/\rho_1\rho_2)^r$,
by Gauss's theorem,
$$
\gamma_n=\frac{(\rho_1)_n(\rho_2)_n}{(aq/\rho_1)_n(aq/\rho_2)_n}
\biggl(\frac{aq}{\rho_1\rho_2}\biggr)^n
\cdot
\prod_{m=1}^\infty\biggl[\frac{(1-aq^m/\rho_1)(1-aq^m/\rho_2)}{(1-aq^m)(1-aq^m/\rho_1\rho_2)}\biggr].
$$
Consequently
\begin{align*}
&
\sum_{n=0}^\infty(\rho_1)_n(\rho_2)_n(aq/\rho_1\rho_2)^n\beta_n
\\ &\qquad\qquad
=\prod_{m=1}^\infty\biggl[\qquad\biggr]
\sum_{n=0}^\infty\frac{(\rho_1)_n(\rho_2)_n}{(aq/\rho_1)_n(aq/\rho_2)_n}
\biggl(\frac{aq}{\rho_1\rho_2}\biggr)^n\alpha_n.
\end{align*}
When $a=1$ or $a=q$ this is equivalent to Rogers's formulae with $u$ \& $v$ in them,%
\footnote{This makes the $uv$ formulae seem almost trivial!}
from which he deduces 19 particular cases giving Fourier series in terms of $A$'s. Actually
$u=q^{\frac12}/\rho_1$, $v=q^{\frac12}/\rho_2$.

Similarly the relations between $\beta_n$ \& $\alpha_n$ gives, with the analogue of Dougall,
relations corresponding to those given by Dougall connecting the $a$'s \& $b$'s.
It is evident  from all this, for example, that any results obtained from Rogers'
formulae E1, E3, F1, F2, E2, E4, F3, F4 (+\,perhaps others) all the 19 $uv$ formulae can be
derived from Watson's transformation directly. Of course one could work out a
few formulae which would generalise all those given by Rogers or obtainable
by the methods of his 1917 paper, but, apart from the results given already, these formulae
appear to me to be anything but elegant. Some have some factors advancing by $\sqrt q$ \& some by $q$,
\& the series are not of any general type. In fact it seems to me that, apart from the
general transformations already given, the method is only useful for obtaining formulae
of the Rogers'--Ramanujan type. They are, at any rate, reasonably simple in appearance.

Of course the simple result at the beginning of this letter has its analogue for integrals. Thus if
\begin{align*}
F(y)&=\int_0^y\varphi(x)f(y-x)g(y+x)\,\d x,
\\ \intertext{and}
G(y)&=\int_y^\infty\psi(x)f(x-y)g(x+y)\,\d x,
\end{align*}
then
$$
\int_0^\infty\varphi(x)G(x)\,\d x
=\int_0^\infty\psi(x)F(x)\,\d x.
$$

One might hope that results could be got for integrals corresponding to those got
for series, but the trouble is to start. So far I have got nowhere.

The result of all this is that I feel that the only thing I am being led to is a
search for more R-R identities, \& probably you have found the most interesting ones that
are new. Of course the method gives $a$-generalisations of them.

Yours sincerely,

$\underline{\text{W.~N.~Bailey}}$

\hbox to\hsize{\kern6mm\hrulefill\kern6mm}
\vskip6mm

\setcounter{section}{6}
\setcounter{equation}{3}

\hbox to\hsize{\hss\hbox{\vtop{\hsize=36mm
\leftline{17 Prince's Road,}
\leftline{Heaton Moor,}
\leftline{Stockport,}
\leftline{Cheshire.}
\smallskip
\leftline{Aug.~1/44}
}}}

\vskip4mm
\noindent
Dear Dyson, \\
\hbox to 24mm{\hfill}%
Hardy has passed on your comments on my paper on ``Identities
of the Rogers--Ramanujan type.'' Many thanks for reading it so carefully \&
for finding the errors. I had checked the formulae a little, but evidently
not enough. The first two formulae were got from
\begin{multline}
1+\sum_{n=1}^\infty\frac{(\rho_1)_n(\rho_2)_n\{ax^3\}_{n-1}}{x_n!(ax)_{2n-1}}\biggl(\frac{ax}{\rho_1\rho_2}\biggr)^n
\\
=\prod_{m=1}^\infty\biggl[\frac{(1-ax^m/\rho_1)(1-ax^m/\rho_2)}{(1-ax^m)(1-ax^m/\rho_1\rho_2)}\biggr]
\\ \times
\biggl[1+\sum_{n=1}^\infty\frac{(-1)^n\{ax^3\}_{n-1}(1-ax^{6n})(\rho_1)_{3n}(\rho_2)_{3n}}{x^3_n!\,(ax/\rho_1)_{3n}(ax/\rho_2)_{3n}}
\\ \times
\frac{a^{4n}x^{\frac32n(3n+1)}}{\rho_1^{3n}\rho_2^{3n}}\biggr],
\end{multline}
where $\{a\}_n=(1-a)(1-ax^3)\dotsb(1-ax)^{3n-3}$.
\setcounter{section}{7} \setcounter{equation}{0}
\\
If we take $\rho_1=-\sqrt a$, $\rho_2=-\sqrt{ax}$, this becomes
\begin{multline*}
1+\sum_{n=1}^\infty\frac{(-\sqrt a)_n(-\sqrt{ax})_n\{ax^3\}_{n-1}x^{\frac12n}}{x_n!(ax)_{2n-1}}
=\prod_{m=1}^\infty\biggl[\frac{(1+\sqrt ax^m)(1+\sqrt ax^{m-\frac12})}{(1-ax^m)(1-x^{m-\frac12})}\biggr]
\\ \times
\biggl[1+\sum_{n=1}^\infty\frac{(-1)^n\{ax^3\}_{n-1}(1-\sqrt ax^{3n})(1+\sqrt a)a^nx^{\frac{9n^2}2}}{x^3_n!}\biggr].
\end{multline*}
I got the first incorrect result by taking $\sqrt a=x^{\frac32}$. I find that I dropped a factor $(-1)^n$, \& the formula
I now get is
\begin{gather}
\label{e7.1}
\sum_{n=0}^\infty\frac{x^6_n!\,x^n}{x_{2n+2}!\,x^2_n!}
=\prod_{n=1}^\infty\frac{(1-x^{18n})(1-x^{18n-3})(1-x^{18n-15})}{(1-x^n)(1-x^{2n-1})}.
\\ \intertext{Similarly by taking $\sqrt a=1$, I got the second formula, viz}
\label{e7.2}
1+2\sum_{n=1}^\infty\frac{x^6_{n-1}!\,x^n}{x_{2n-1}!\,x^2_n!}
=\prod_{n=1}^\infty\frac{(1+x^n)(1-x^{9n})}{(1-x^n)(1+x^{9n})}.
\end{gather}
You say this is correct up to the term $x^6$. I must confess that I cannot find
anything wrong in the working.

The third formula was certainly wrong, \& should have been the same as the first of
the five you sent. I will incorporate these formulae in the paper, but I should be very glad
if you could say whether you agree with \eqref{e7.1} \& \eqref{e7.2} now. I find these things
very tedious to check to any extent, though I thought I was fairly safe.

Again, many thanks for all the care you have taken \& for the new results.

Yours sincerely,

$\underline{\text{W.~N.~Bailey}}$

\hbox to\hsize{\kern6mm\hrulefill\kern6mm}
\vskip6mm

\hbox to\hsize{\hss\hbox{\vtop{\hsize=36mm
\leftline{8 Langton Avenue,}
\leftline{Whetstone,}
\leftline{London, N.~20.}
\smallskip
\leftline{Oct.~8/46}
}}}

\vskip4mm
\noindent
Dear Dyson, \\
\hbox to 24mm{\hfill}%
I was interested to hear from you again \& that you are back to Cambridge.
I am now in London (at Bedford College) \& have been for the past two years.

I have had two papers ready for P.L.M.S.\ for about 3 years or more, so they ought
to be published in another year or two!
Actually the L.M.S.\ have done all they can to speed up publication, but first of all
shortage of paper \& then shortage of labour have been too much for them.

After coming here I had rather a strenuous time getting used to the ways of London
University, finding a house, \& so on, so I didn't make much progress with the work I was doing.
Lately, however, I have sent a paper to the Quarterly Journal in which I give the basic
analogues of 6.6(3) in my tract, \& of 6.8(3) \& 7.6(2). The first of these is what
Dougall's theorem becomes when the series does not terminate. The other two are
the relations connecting 4 ${}_9F_8$'s. I found the analogue of 6.6(3) was merely
a particular case of a formula given in a Q.J.\ paper in 1936 (Series of hyp.\ type
infinite in both dir$^{\text{ns}}$, Q.J.\ 7 (105) first formula in \S\,5).

I got the analogue of 6.8(3) by transforming the argument in \S\S\,6.7 \& 6.8 of my tract
into series by considering poles on the right of the contours. Then I did the corresponding
work for basic series. The idea was simple enough, but the details nearly broke my heart.

In P.L.M.S.\ 42 (1934) 410--421, Whipple gave (or rather showed how to find) a connection
between 4 ${}_9F_8$'s (well-poised) when there was no restriction on the sums
of numerator \& den.\ parameters. This leaves the obvious problem of finding the
corresponding result for basic series, but I hadn't the pluck to start that. Whipple's
proof is very short \& depends on a contour integral, but I don't see how one could
adapt this method to basic series, unless one worked out a good deal about
integrals generalising integrals of Barnes' type.

These was another thought I had that seemed to hold promise at one time,
but I never got anything out of it. In the papers you saw I gave a general
theorem on series which has the integral analogue: If
\begin{align*}
F(y)&=\int_0^y\varphi(x)f(y-x)g(y+x)\,\d x
\\ \intertext{\&}
G(y)&=\int_y^\infty\psi(x)f(x-y)g(x+y)\,\d x,
\end{align*}
then
$$
\int_0^\infty\varphi(x)G(x)\,\d x
=\int_0^\infty\psi(x)F(x)\,\d x.
$$
With so much being derivable from the series theorem I thought that there
might be possibilities from the integral theorem, but I did not succeed in
getting anything interesting. Still, there may be something in it.

I didn't try to get any more identities of the Rogers--Ramanujan type.
We got a good many between us 3 years ago!

Yours sincerely,

$\underline{\text{W.~N.~Bailey}}$

\hbox to\hsize{\kern6mm\hrulefill\kern6mm}


\end{document}